\documentclass{amsart}
\usepackage{amscd,amsmath,amssymb}
\usepackage[mathscr]{eucal}
\usepackage{graphicx}
\newtheorem{theorem}{Theorem}[section]

\newtheorem{corollary}[theorem]{Corollary}
\newtheorem{lemma}[theorem]{Lemma}
\newtheorem{proposition}[theorem]{Proposition}
\theoremstyle{definition}
\newtheorem{definition}[theorem]{Definition}
\theoremstyle{remark}
\newtheorem{remark}[theorem]{Remark}

\numberwithin{equation}{section}

\newcommand{\eps}{\varepsilon}

\begin{document}
\title{On the extrinsic geometry of contact structures.}
\author{Vladimir Krouglov.}
\maketitle
\begin{abstract}
In the paper we prove, that extrinsic curvature does not impose restrictions on the topology of a contact structure, except the obvious ones.
\end{abstract}

\section{Introduction.}
One of the natural characteristics of the plane distribution on a Riemannian manifold is its second fundamental form. It is known, a foliation on a closed three manifold is taut if the trace of the second fundamental form (i.e. the mean curvature) vanishes with respect to some Riemannian metric. It can be shown \cite{YC}, that vanishing of the second fundamental form of a distribution with respect to some metric yields topological restrictions on the ambient space. In the recent work by Patrick Massot \cite{PM} it was established that all totally geodesic contact structures on three manifolds are tight. In contrast with foliations, every contact structure is a minimal distribution. Therefore it is an interesting question how one can relax the condition of a contact structure being totally geodesic to keep the topological restrictions on a contact structure.

In the present paper we study contact structures on 3-manifolds which have the restrictions on the determinant of the second fundamental form (i.e. the extrinsic curvature). There are three natural classes of plane distributions with the restrictions on the extrinsic curvature:
\begin{definition}
The distribution of planes $\xi$ on a three manifold is called:
\begin{enumerate}
\item{parabolic, if there is a Riemannian metric on $M$ such that the extrinsic curvature of $\xi$ is equal to zero.}
\item{strong saddle (or hyperbolic), if there is a Riemannian metric on $M$ such that the exrinsic curvature of $\xi$ is strictly less than zero.}
\item{elliptic, if there is a Riemannian metric on $M$ such that the extrinsic curvature of $\xi$ is everywhere positive.}
\end{enumerate}
\end{definition}
The main result of this paper is the following:
\begin{theorem}
Let $\xi$ be a transversally orientable contact structure on a closed orientable three manifold $M$. Then,
\begin{enumerate}
\item{$\xi$ is parabolic.}
\item{$\xi$ is hyperbolic if and only if its Euler class vanishes.}
\item{There are no elliptic plane distributions on closed three manifolds.}
\end{enumerate}
\end{theorem}
\begin{remark}
In \cite{Kru1} and \cite{Kru2} we studied the problem of existence of foliations with the restrictions on the extrinsic curvature of the leaves. The main result in \cite{Kru2} is obtained in the present paper from some different point of view.
\end{remark}
This paper is organized as follows. In Section $2$ we recall basic definitions and results in the geometry of plane distributions. Notions and results from the contact topology and Giroux correspondence between contact structures and open book decompositions are reviewed in Section $3$. Section $4$ is devoted to the proof of the fact, that every foliation by the fibers of fibration over the circle is parabolic. We also define a special parabolic foliation associated with the open book decomposition of $M$. In Section $5$ we give the proof of the main theorem.
\section{Basic Definitions and Notations.}
Throughout this paper $M$ will be a closed orientable 3-manifold. A distribution on $M$ is a two dimensional subbundle of the tangent bundle of $M$. That is, at each point $p$ in $M$ there is a plane $\xi_p$ in the tangent space $T_p M$. A distribution is called integrable, if there is a foliation on $M$ which is tangent to it. The following  theorem of Frobenius gives necessary and sufficient conditions for $\xi$ to be integrable.
\begin{theorem}
Let $\xi$ be a distribution on $M$. Then $\xi$ is integrable if and only if for any two sections $S$ and $T$ of $\xi$ its Lie bracket belongs to $\xi$.
\end{theorem}
\begin{definition}
A distribution $\xi$ is called a contact structure if for any linearly independent sections $S$ and $T$ of $\xi$ and for any $p \in M$ the Lie bracket $[S, T]$ at $p$ does not belong to $\xi_p$.
\end{definition}
A distribution $\xi$ is called transversally oriented if there is a globally defined $1$-form $\alpha$, such that $\xi = Ker(\alpha)$. This is equivalent to say that there exists a globally defined vector field $n$ which is transverse to $\xi$. It is an easy consequence of Frobenius Theorem that $\xi$ is a contact structure if and only if
$$
\alpha \wedge d\alpha \ne 0
$$
The Euler class $e(\xi) \in H^2(M,\mathbb{Z})$ of a plane distribution is the Euler class of the bundle $\xi \to M$.
It is known that if $\xi$ is a 2-dimensional plane distribution on $M$ with vanishing Euler class then $\xi$ is trivial.

Assume that $M$ is a Riemannian manifold with the metric $\langle \cdot, \cdot \rangle$ and the Levi-Civita connection $\nabla$. Let $n$ be a local unit vector field orthogonal to $\xi$. We are now going to define the second fundamental form of $\xi$. The definition is due to Reinhart \cite{Re}.
\begin{definition}
The second fundamental form of $\xi$ is a symmetric bilinear form, which is defined in the following way:
$$
B(S, T) = \frac{1}{2}\langle \nabla_S T + \nabla_T S, n \rangle
$$
for all sections $S$ and $T$ of $\xi$.
\end{definition}
\begin{remark} If $\xi$ is integrable, then $B$ restricted to the leaf of $\xi$ agrees with the second fundamental form of the leaf.
\end{remark}
Assume that $S$ and $T$ are two linearly independent sections of $\xi$.
\begin{definition}
A mean curvature function $H$ of a plane distribution is simply  a trace of the second fundamental form. If $S$ and $T$ are orthonormal, it may be written as
$$
H = \langle \nabla_S S, n\rangle + \langle \nabla_T T, n\rangle
$$
\end{definition}
\begin{definition}We call the function
$$
K_e(\xi) = \frac{B(S, S)B(T, T) - B(S, T)^2}{\langle S,S\rangle \langle T,T\rangle - \langle S, T\rangle^2}
$$ an extrinsic curvature of $\xi$.
\end{definition}
It is easy to verify that $K_e(\xi)$ depends only on $\xi$, not on the actual choice of $S$, $T$ and $n$.

There are several classes of plane distributions on 3-manifolds in depending on the signs of the extrinsic an mean curvature of $\xi$.
\begin{definition}A plane distribution on a three manifold is called
\begin{enumerate}
\item{totally geodesic, if there is Riemannian metric on $M$ such that $B \equiv 0$.}
\item{minimal, if there is a Riemannian metric on $M$ such that $H = 0$}
\item{ parabolic, if there is a Riemmanian metric on $M$, such that the extrinsic curvature $K_e(\xi) = 0$.}
\item{strong saddle (or hyperbolic), if there is a Riemmanian metric on $M$ such that the extrinsic curvature $K_e(\xi) < 0$. }
\item{elliptic, if there is a Riemmanian metric on $M$ such that the extrinsic curvature $K_e(\xi) > 0$.}
\end{enumerate}
\end{definition}
\section{Open book decompositions and contact structures.}
Consider an oriented link $L$ in an oriented three manifold $M$. Assume that a complement $M \backslash L$ fibers over the circle with the projection map $\pi : M\backslash L \to S^1$ and that $\pi^{-1}(t)= \Sigma^2_t$ is an interior of a compact surface bounded by $L$. Pair $(\pi, L)$ is called an open book decomposition of $M$.

Two open book decompositions $(L, \pi)$ and $(L', \pi')$ are called isomorphic if there is a diffeomorphism $f: M \to M$ such that $\pi' \circ f = \pi$.

There is an alternative description of open book decompositions through the mapping cylinders. Assume that $\Sigma^2$ is a compact surface (with boundary) and consider a mapping cylinder
$$\Sigma^2 \times_\phi S^1 = \Sigma^2 \times [0, 1] \slash  (x, 0) \sim (\phi x, 1),$$
where $\phi$ is some diffeomorphism of $\Sigma^2$ which is an identity in the neighborhood of $\partial \Sigma^2$. Since $\phi$ is an identity map in the neighborhood of the boundary the boundary of $\Sigma^2 \times_\phi S^1$ is a union of $r$ tori where $r$ is a number of connected components of the boundary $\partial \Sigma^2$. We may now glue $r$ solid tori $D^2 \times S^1$ to $\Sigma^2 \times_\phi S^1$ in such a way that $\partial D^2$ is glued to $S^1 = [0, 1] \slash \sim$ and $S^1$ factor in the solid torus $D^2 \times S^1$ corresponds to the boundary component of $\Sigma^2$. As a result we obtain closed manifold $M = (\Sigma^2 \times_\phi S^1)\cup (\bigcup_r D^2 \times S^1)$. This manifold has a canonical presentation as an open book decomposition such that $L$ is a union of $r$ core curves of $D^2 \times S^1$ that we glued to $\Sigma^2 \times_\phi S^1$ to obtain $M$.

\begin{definition}
We say that contact structure $\xi$ is supported by the open book decomposition $(L, \pi)$ of $M$ if there is a one-form $\alpha$ associated with  $\xi$ such that $\alpha(L) > 0$ and $d\alpha|_\Sigma^2$ on each page $\Sigma^2 = \pi^{-1}(t)$.
\end{definition}
In \cite{ThW}, Thurston and Winkelhemper have shown that each open book decomposition on a closed three manifold supports some contact structure. It is surprising that converse also holds.
\begin{theorem}(Giroux, \cite{Gir})
Every contact structure $\xi$ on a closed orientable three manifold $M$ is supported by some open book decomposition.
\end{theorem}
Using this result Etnyre in \cite{Et} proved the following result:
\begin{theorem}
Every contact structure on a closed orientable three manifold is a $C^\infty$-deformation of $C^\infty$-foliation.
\end{theorem}
Consider the diffeomorphism $\phi$ of a compact surface $\Sigma^2$ that is an identity in the neighborhood of $\partial \Sigma^2$ (or an arbitrary diffeomorhpism if $\Sigma^2$ is closed). Throughout the paper $\mathcal{F}_\phi$ will denote the foliation of the mapping cylinder $\Sigma^2 \times_\phi S^1$ by the surfaces $\Sigma^2 \times \{t\}$ for all $t \in S^1$.

\section{Parabolic foliations.}
In \cite{Kru2} we showed that each closed orientable three manifold admits a parabolic foliation (i.e. the foliation by parabolic surfaces with respect to some metric). We established this result using Dehn surgeries along knots in $S^3$. In the present section we will show how one can obtain this result using open book decompositions. We will also show that every foliation by the fibers of fibration over the circle is parabolic.
\subsection{Local models of parabolic foliations.}
The following lemma was proved in â \cite{Kru2}. This lemma allows one to glue two parabolic foliations together along the common boundary leaf preserving parabolicity of glued foliation.
\begin{lemma}\label{lem:2} Let $\Sigma^2$ be a compact parallelizable surface (possibly with the boundary). Consider two metrics $G$ and $H$ on $\Sigma^2$ that coincide in some neighborhood of the boundary. Let $\mathcal{F}$ be a foliation of $M = \Sigma^2 \times [0, 1]$ by the surfaces $\Sigma^2 \times \{t\}$. Then, there is a metric $g$ on $M$ such that
\begin{enumerate}
\item{In some tubular neighborhood of $\Sigma^2 \times \{0\}$, $g = dt^2 + G(p)$ for all $p \in \Sigma^2$.}
\item{in some tubular neighborhood of $\Sigma^2 \times \{1\}$, $g = dt^2 + H(p)$ for all $p \in \Sigma^2$.}
\item{$\mathcal{F}$ is a parabolic foliation on $\Sigma^2 \times [0, 1]$ with respect to $g$}.
\item{There is a neighborhood $U$ of the boundary $\partial \Sigma^2$ such that for all  $t \in [0,1]$, $$g(p, t)|_{U \times \{t\}} = G(p)$$}
\end{enumerate}
\end{lemma}
The next lemma shows that a Reeb foliation inside a solid torus is parabolic. The result is due to Bolotov \cite{Bol}.
\begin{lemma}\label{lem:1}
There is a foliation $\mathcal{F}$ and a metric $g$ on $D^2 \times S^1$ such that:
\begin{enumerate}
\item{$\mathcal{F}$ is parabolic with respect to $g$.}
\item{The foliation $\mathcal{F}|_{D^2(\frac{1}{3}) \times S^1}$ is a foliation by the totally geodesic disks $D^2(\frac{1}{3}) \times \{t\}$ and the foliation $\mathcal{F}|_{([\frac{2}{3}, 1] \times S^1) \times S^1}$ is a foliation by the totally geodesic tori $\{ r \} \times S^1 \times S^1$.}
\end{enumerate}
\end{lemma}
\emph{Proof:} Consider the solid torus $D^2 \times S^1$ with the following coordinates on it
$$
D^2 \times S^1 = \{ ((r, \phi), t) : r \in [0, 1], \phi, t \in [0, 2\pi) \}
$$
Define  the one-from $\alpha$ on $D^2 \times S^1$ as:
$$
 \alpha = f(r)dr + (1 - f(r)) dt
$$
where $f(r)$ is such a smooth function on $[0,1]$ that
$$
 f(r) = \left \{ \begin{array}{l} 0, \ r \in [0, \frac{1}{3}] \\
                                  \mbox{is a strictly increasing function when} \   r \in (\frac{1}{3}, \frac{2}{3}] \\
                                  1, \   r \in (\frac{2}{3}, 1]
  \end{array} \right.
$$
This form defines a `thick' Reeb foliation $\mathcal{F}$ on $D^2 \times S^1$ (that is, there is a subset $N$ such that $\mathcal{F}|_{N}$ is a Reeb foliation and $\mathcal{F}|_{D^2 \times S^1 \backslash N}$ is diffeomorphic to a product foliation by tori).

Assume that in coordinates $(r, \phi, t)$ the matrix of $g$ has a form:
$$
 g = \left ( \begin{array}{ccc} 1 & 0 & 0 \\
 0 & G(r) & 0 \\
 0 & 0 & 1 \end{array} \right )
$$
In order to calculate the second fundamental form of $\mathcal{F}$ consider the following sections: $X = \frac{\partial}{\partial \phi}, Y =(1 - f(r))\frac{\partial}{\partial r} - f(r) \frac{\partial}{\partial t}$ of the tangent bundle $T\mathcal{F}$. Let $n = f(r) \frac{\partial}{\partial r} + (1 - f(r))\frac{\partial}{\partial t}$ be a normal vector field.

By the straightforward calculation we obtain that the matrix of the second fundamental form is equal to:
$$
\frac{1}{2f(r)^2 - 2f(r) +1}\left ( \begin{array} {cc} - f \frac{\partial G}{\partial r} & 0 \\
                           0 & - (1-f) \frac{\partial f}{\partial r} \end{array} \right )
$$
It is obvious that since $f = 0$ on $[0, \frac{1}{3})$ the foliation is by totally geodesic disks for every choice of $G=G(r)$. Define $G = G(r)$ in the following way:
$$
G = \left \{ \begin{array} {l} r^2, \mbox{when $r \in [0, \frac{1}{4})$} \\
                               \mbox{strictly increasing, when $r \in [\frac{1}{4}, \frac{1}{3})$} \\
                               \mbox{$1$, when $r \in [\frac{1}{3}, 1]$}   \end{array} \ \right.
$$
For this choice of $G$, the metric $g$ is regular in the neighborhood of the core curve $r = 0$ and satisfies conditions of the lemma.

\subsection{Parabolic fibrations.}
Let $(M^3, \pi)$ be a fibration over the circle with a closed fiber $\Sigma^2$ and let $\mathcal{F}_\phi$ be the foliation by the fibers of $\pi$. We have the following
\begin{proposition}\label{prop:1}
The foliation $\mathcal{F}_\phi$ is parabolic.
\end{proposition}
\emph{Proof:} Consider the presentation of $M$ as a mapping cylinder of some $\phi: \Sigma^2 \to \Sigma^2$. We may assume that $\phi$ is a composition of Dehn twists $\alpha_i$ along some closed curves $\gamma_i$ in $\Sigma^2$ and that support of each $\alpha_i$ is contained in some annulus $c_i$ for each $\ i$.

First, consider the case when $\phi$ is a Dehn twist by itself. Pick an arbitrary metric $G$ on $\Sigma^2$. On a manifold $N = \Sigma^2 \times [0, 1]$ consider a direct product metric $G + dt^2$. It is obvious that a foliation of $N$ by the surfaces $\Sigma^2 \times \{t\}$ is parabolic with respect to this metric.

Denote by $H$ the pullback metric $\phi^\ast G$. Since the support of $\phi$ is contained in some annulus $c$, from Lemma $4.1$ we may obtain that on $c\times [0, 1]$ there is a metric $g$ such that $g|_{c\times \{0\}} = dt^2 + G|_c$, $g|_{c\times \{1\}} = dt^2 + H|_c$ and foliation by $c\times \{t\}$ is parabolic. From condition $4)$ of Lemma $4.1$ this metric is glued smoothly with the direct product metric $G + dt^2$ on the complement $(\Sigma^2 \backslash c) \times [0, 1]$. Consequently, on $L = \Sigma^2 \times [0, 1]$ there is a metric (which will also be denoted by $g$) that a foliation by surfaces $\Sigma^2 \times \{t\}$ is parabolic with respect to it. We are left to consider the union
$$
M = N \cup_\phi L
$$
and glue the remaining boundary components of $M$ by the identity map. It is clear that since $g$ is a direct product metric in the one-sided neighborhoods of the boundary of $L$ (and $\frac{\partial}{\partial t}$ is a unit normal vector field), it is glued correctly with the direct product metric on $N$ and therefore define a smooth metric on $M$. The foliation $\mathcal{F}_\phi$ is parabolic with respect to the introduced metric.

Assume now that $\phi = \alpha_1 \circ \alpha_2 \circ \ldots \circ \alpha_n$. Pick an arbitrary metric $G$ on $\Sigma^2$ and for each  $i$ consider the manifold $N_i = \Sigma^2 \times [0, 1]$ with the direct product metric $dt^2 + G$. For each $N_i$ consider the corresponding $L_i$ that is obtained as on the previous step. In the (one-sided) neighborhood of the boundary components of $\partial(N_i \cup_{\alpha_i} L_i)$ the metric is a direct product metric $dt^2 + G$. We are left to consider the union
$$
M = (N_1 \cup_{\alpha_1} L_1) \cup_{id} (N_2 \cup_{\alpha_2} L_2) \cup_{id} \ldots \cup_{id}(N_n \cup_{\alpha_n} L_n)
$$
and glue the remaining boundary components by the identity map. This finishes the proof of proposition.
\begin{corollary}\label{cor:1}
Let $\Sigma^2$ be a compact orientable surface with boundary and assume that $\phi$ is a diffeomorphism of $\Sigma^2$ which is an identity map in the neighborhood of the boundary $\partial \Sigma^2$. On the mapping cylinder $\Sigma^2 \times_\phi S^1$ the foliation by surfaces $\Sigma^2 \times \{t\}$ is parabolic.
\end{corollary}
\emph{Proof:} Same as Proposition \ref{prop:1}
\subsection{Parabolic foliations associated with the open book decomposition.}
Assume that $(L, \pi)$ is an open book decomposition of the closed oriented three manifold $M$. In this section we will define a special parabolic foliation on $M$ associated with this open book decomposition. We will define this foliation replacing the neighborhoods of the binding by a foliation which was defined in Lemma \ref{lem:1} and twisting the pages of the open book around the boundary leaf.

Denote by $N$ the tubular neighborhood of one component of the binding (construction for other components of the binding would be the same). On $N$ define the coordinates $(r, \phi, t)$ in such a way that pages of the open book correspond to the constant $\phi$-annuli and the component of the binding corresponds to a curve $r = 0$. Assume that
$$
N = \{(r, \phi, t) : r \le 1 + 2\eps \}
$$
for some small fixed $\eps$.

Define a foliation and the metric  on $N_1 = \{(r, \phi, t) : r\le 1\}$ as in Lemma \ref{lem:1}.
Consider the following function $f(r)$ on $[1, 1+2\eps]$:
\begin{enumerate}
\item{$f(r) = 0$, for all $r \in [1, 1+ \frac{\eps}{2}]$}
\item{$f(r) = 1$, for all $r \in [1+\eps, 1+2\eps]$}
\item{$f$ is a strictly increasing function on $[1+\frac{\eps}{2}, 1+\eps]$}
\end{enumerate}
On $[1, 1+2\eps]\times T^2$ define a foliation using a one-form:
$$
\alpha = f(r) d\phi + (1 - f(r))dr
$$
This foliation when restricted to $[1 + \eps, 1+ 2\eps] \times T^2$ is a foliation by the constant $\phi$-annuli and is therefore smoothly glued to a foliation of $M \backslash N$ by the fibers of $\pi$.

In the neighborhood of the boundary $\partial N_1$ the metric $g = dr^2 + d\phi^2 + dt^2$. Extend this metric to $N_{1+3\eps}$. The foliation in $[1, 1+2\eps] \times T^2$ is parabolic with respect to this metric, since in each point of $[1, 1+ 2\eps] \times T^2$ the plane tangent to the leaf contains principal geodesic direction $\frac{\partial}{\partial t}$. Moreover, this metric induces some metric in the neighborhood of one component of the  boundary of $\Sigma^2 \times_\phi S^1$. Making such extensions for each component of the binding we will obtain some metric that is defined in the neighborhood of the boundary $\partial (\Sigma^2 \times_\phi S^1)$. Note, that this metric is direct product metric in this neighborhood of the boundary. It induces some metric in the neighborhood of $\partial \Sigma^2 \times \{0\}$. Arbitrarily extend it to the whole leaf $\Sigma^2\times \{0\}$. Using the same construction as in Lemma \ref{lem:2} we may define a metric on a mapping cylinder in such a way that a foliation by the fibers of $\pi$ is parabolic. Therefore we defined a metric on $\Sigma^2 \times_\phi S^1 \cup (\bigcup_i^r N_{1+2\eps}(i))$, where $r$ is a number of components of the binding in the open book decomposition and the foliation $\mathcal{F}_\phi$ which is parabolic with respect to it.
\section{Extrinsic geometry of contact structures.}
In this section we will give the proof of Theorem $1.2$.\\
\emph{Proof of the Theorem:}
\begin{lemma}\label{lem:3}
The function of the mean curvature of any transversally orientable distribution $\xi$ on a closed orientable three manifold is a divergence.
\end{lemma}
\emph{Proof:} Let $(X, Y)$ be two (local) orthonormal sections of $\xi$ and $n$ be a unit normal vector field. From the definition of the mean curvature function we have
$$
H = tr (B) = \langle \nabla_X X, n \rangle + \langle \nabla_Y Y, n \rangle
$$
Using the fact that $n$ is orthogonal to $\xi$ we obtain that
$$
H = - \langle \nabla_X n, X \rangle - \langle \nabla_Y n, Y \rangle = div(-n)
$$
since $\langle \nabla_n n, n \rangle = 0$.
\begin{corollary}\label{cor:2}
On a closed orientable three manifold there are no transversally orientable elliptic distributions.
\end{corollary}
\emph{Proof:} From the Stokes theorem
$$
\int_M H = \int_M div(-n) = 0
$$
If a distribution is elliptic, then $K_e > 0$. Therefore principal curvature functions have to be either both positive or negative. Consequently, $\int_M H \ne 0$. A contradiction.
\begin{remark}
Lemma \ref{lem:3} and the Corollary \ref{cor:2} are known for a long time. However, author does not know when this fact was used for the first time.
\end{remark}
\begin{lemma}
Transversally orientabe contact structure $\xi$ on a closed orientable three manifold is strong saddle if and only if $e(\xi)=0$
\end{lemma}
\emph{Proof:} See Corollary $3.6$ in \cite{Kru3} for details.

We are left to prove that every transversally orientable contact structure on $M$ is parabolic. The following lemma allows to carry the information about the extrinsic geometry from one distribution to another.
\begin{lemma}\cite{Bol2}
Assume that on a Riemannian manifold $M$ given two distributions $\xi$ and $\eta$. Denote by $B_\xi$ the second fundamental form of $\xi$. Assume that a normal vector field to $\xi$ is transverse to $\eta$. For all $X \in \xi$, denote by $PX$ an orthogonal projection of vector $X$ on $\eta$. Then, there is a metric on $M$ such that a second fundamental form of $\eta$ satisfies
$$
B_\eta(PX, PY) = B_\xi(X, Y), \ \mbox{for all $X$,  $Y \in \xi$.}
$$
In particular, if $\xi$ is parabolic, so is $\eta$.
\end{lemma}
\begin{lemma}
Every transversally orientable contact structure $\xi$ on a closed orientable three manifold $M$ is parabolic.
\end{lemma}
\emph{Proof:} We will follow the construction in \cite{Et} to define a family of contact one-forms that approximate a parabolic foliation associated with the open book decomposition.

For $\xi$ consider the corresponding open book decomposition $(\Sigma^2, \phi)$. It is obvious that foliation associated with this open book decomposition differs from the one defined in \cite{Et} (cf. Proof of Theorem 1 there) everywhere except some $\delta$-neighborhood of the leaf $\{r=1\}$. It is easy to verify that a family of one-forms, defined in \cite{Et} is an approximation of our foliation. So, there is a family of one-forms $\alpha_t, t>0$ which is a deformation of $\mathcal{F}_\phi$ and has support $(\Sigma^2, \phi)$. In particular, $\xi_t = Ker(\alpha_t)$ are isotopic to $\xi$.

Let $g$ denotes a Riemannian metric on $M$ such that $\mathcal{F}_\phi$ is parabolic. It is obvious that we may find such $t>0$ that $\xi_t$ will be transverse to a normal vector field to $\mathcal{F}_\phi$. Using Lemma $5.5$ we see that there is a Riemannian metric on $M$ such that $\xi_t$ is parabolic (for $t$ sufficiently small). But $\xi_t$ is contactomorphic to $\xi$ and therefore $\xi$ is a parabolic contact structure with respect to a pullback metric. This finishes the proof of the theorem.

\end{document}